\documentclass[a4paper,reqno,9pt]{amsart}

\usepackage{amsmath}
\usepackage{amsfonts}
\usepackage{palatino}
\usepackage{a4wide}
\usepackage{amssymb}
\usepackage{amsthm}
\usepackage{epsfig}

\usepackage[usenames]{color}
\definecolor{citegreen}{rgb}{0,0.6,0}
\definecolor{refred}{rgb}{0.8,0,0}

\theoremstyle{plain}

\newtheorem{teo}{Theorem}[section]

\newtheorem{prop}[teo]{Proposition}

\newtheorem{ackn}{Acknowledgments\!}

\theoremstyle{definition}

\theoremstyle{remark}
\newtheorem{rem}[teo]{Remark}

\numberwithin{equation}{section}

\def\R{{{\mathbb R}}}
\def\SS{{{\mathbb S}}}

\def\RRR{{\mathrm R}}

\def\HHH{{\mathrm H}}

\def\Ric{{\mathrm {Ric}}}

\def\QQQ{\mathrm {Q}}
\def\KKK{\mathrm {K}}
\def\kkk{\mathrm {k}}

\begin{document}

\title[Flow by Mean Curvature inside a Moving Ambient
Space]{Flow by Mean Curvature inside a Moving Ambient
Space}

\author[Annibale Magni]{Annibale Magni}
\address[Annibale Magni]{Mathematisches Institut, Abt. f\"ur Reine
  Mathematik, Albert--Ludwigs--Universit\"at, Eckerstrasse 1, D-79104, Freiburg im Breisgau, Germany}
\email[A. Magni]{annibale.magni@math.uni-freiburg.de}
\author[Carlo Mantegazza]{Carlo Mantegazza}
\address[Carlo Mantegazza]{Scuola Normale
  Superiore di Pisa, Piazza dei Cavalieri 7, Pisa, Italy, 56126}
\email[C. Mantegazza]{c.mantegazza@sns.it}
\author[Efstratios Tsatis]{Efstratios Tsatis}
\address[Efstratios Tsatis]{Department of Mathematics, University of Patras, Patras 26500, Greece}
\email[Efstratios Tsatis]{etsatis@upatras.gr}

\date{\today}

\begin{abstract} We show some computations related to the
  motion by mean curvature flow of a submanifold inside an ambient
  Riemannian manifold evolving by Ricci or backward Ricci flow. 
  Special emphasis is given to the possible generalization of Huisken's monotonicity
  formula and its connection with the validity of some Li--Yau--Hamilton differential 
Harnack--type inequalities in a moving Riemannian manifold.
\end{abstract}

\maketitle

%\tableofcontents

\section{Introduction}

In this paper we present some computations concerning the 
mean curvature flow of a submanifold inside a moving Riemannian
manifold. We are particularly interested in finding analogues of
Huisken's monotonicity formula. 
We will see that in some special situations, notably, when 
the ambient is a {\em gradient Ricci soliton}, such a monotonicity 
actually holds (see Section~\ref{solsec}). 
We will analyze in detail the cases when the ambient Riemannian
manifold evolves by Ricci or backward Ricci flow and
we will discuss the connection between the monotonicity of Huisken's
integral and the validity of some Li--Yau--Hamilton differential 
Harnack--type inequalities in the moving manifold.

Some of the computations here are mentioned by Ni~\cite{leini4}. 
A very closely related paper is the one by Lott~\cite{lott2}. 
Moreover, the work of Ecker~\cite{eck3} and the discussion 
in Section~3.10, Chapter~11 of the book by Chow,
Lu and Ni~\cite{chowluni} also deal with the subject of coupling the Ricci
flow with the mean curvature flow.

We recall the fundamental Huisken's monotonicity 
formula (see ~\cite{huisk3}) for the mean curvature flow (from now on
MCF) in the Euclidean space.\\
Let us assume that we have a smooth, compact, $n$--dimensional 
submanifold $N$ immersed in $\R^m$ evolving by the MCF. That is, the flow is described by a smooth map 
$\varphi:N\times[0,T)\to\R^m$ with 
$$
\partial_t\varphi(p,t)=\HHH(p,t)
$$
where the map $\varphi_t=\varphi(\cdot,t):N\to\R^m$ is an immersion for every
$t\in[0,T)$. Here $\HHH(p,t)$ is the vector valued mean
curvature of the submanifold at time $t$ and point $p$.\\
The immersion $\varphi_t$ induces (by pull--back of the standard
scalar product of $\R^m$) a metric $h_t$ on $N$ at every
time $t$, turning $(N,h_t)$ into a Riemannian manifold with a
canonically associated Riemannian volume measure $\mu_t$.\\
We then consider the {\em backward} heat kernel $\rho_{x_0,T}(x,t)$ on
$\R^m$ centered at some point $x_0\in\R^m$ and with maximal time $T>0$, that is,
$$
\rho_{x_0,{T}}(x,t)=\frac{e^{-\frac{\vert
      x-x_0\vert^2}{4({T}-t)}}}{[4\pi({T}-t)]^{m/2}}\,.
$$

\begin{teo}[Huisken's Monotonicity Formula~\cite{huisk3}]
For every $x_0\in\R^m$ and ${T}>0$, there holds
\begin{align}
 \frac{d\,}{dt}\,\Bigl\{[4\pi(T-t)]^{\frac{m-n}{2}}\int_N \rho_{x_0,{T}}\,d\mu_t\Bigr\} 
=&\, \frac{d\,}{dt}\,\int_N \frac{e^{-\frac{\vert
      x-x_0\vert^2}{4({T}-t)}}}{[4\pi({T}-t)]^{n/2}}\,d\mu_t\label{huiskform}\\
=&\,-\int_N 
\left\vert\HHH+\frac{(x-x_0)^\perp}{2({T}-t)}\right\vert^2\,\frac{e^{-\frac{\vert
      x-x_0\vert^2}{4({T}-t)}}}{[4\pi({T}-t)]^{n/2}}\,d\mu_t\nonumber\\
=&\,-[4\pi(T-t)]^{\frac{m-n}{2}}\int_N \left\vert\HHH-\nabla^\perp\log{\rho_{x_0,{T}}}\right\vert^2\, \rho_{x_0,{T}}\,d\mu_t\nonumber
\end{align}
in the time interval $[0,T)$, where $\nabla^\perp$ denotes the projection on the normal space to $N$ 
of the gradient in $\R^m$ of a function.\\
Hence, the integral $\int_N \frac{e^{-\frac{\vert
      x-x_0\vert^2}{4({T}-t)}}}{[4\pi({T}-t)]^{n/2}}\,d\mu_t$ is
nonincreasing during the flow in $[0,T)$.
\end{teo}

Following Hamilton~\cite{hamilton8}, we can consider, more generally, the flow by mean curvature 
of an $n$--dimensional, smooth, compact submanifold $N$ of a Riemannian manifold $(M,g)$ in a  time interval $[0,T)$ and a positive solution $u:M\times[0,T)\to\R$ of the backward heat
equation $u_t=-\Delta^M u$ in the "ambient" space.\\
Making use of the formula
$$
\Delta^Nu=\Delta^M u-g^{\alpha\beta}\nabla^2_{\alpha\beta}u
+\langle\nabla^M u\,\vert\,\HHH\rangle\,,
$$
where we denoted the "normal" indices with Greek letters (this means
that the intermediate term on the right hand side is the ``trace'' of
the 2--form $\nabla^2u$, restricted only to the normal
space to $N$), we compute
\begin{align*}
  \frac{d\,}{dt}\,\Bigl\{[4\pi&\,({T}-t)]^{\frac{m-n}{2}}\int_N u\,d\mu_t\,\Bigr\}\\
  =&\,-[4\pi({T}-t)]^{\frac{m-n}{2}}\int_N \Bigl(-u_t
  +\vert\HHH\vert^2 u -\langle \nabla^M
  u\,\vert\,\HHH\rangle+\frac{(m-n)}{2({T}-t)}u\Bigr)\,d\mu_t\\
=&\,-[4\pi({T}-t)]^{\frac{m-n}{2}}\int_N \Bigl(\Delta^Mu
+\vert\HHH\vert^2 u -\langle \nabla^M
u\,\vert\,\HHH\rangle+\frac{(m-n)}{2({T}-t)}u\Bigr)\,d\mu_t\\
=&\,-[4\pi({T}-t)]^{\frac{m-n}{2}}\int_N \Bigl(\Delta^Nu+g^{\alpha\beta}\nabla^2_{\alpha\beta}u
+\vert\HHH\vert^2 u -2\langle \nabla^M u\,\vert\,\HHH\rangle+\frac{(m-n)}{2({T}-t)}u\Bigr)\,d\mu_t\,.
\end{align*}
As the integral of $\Delta^Nu$ is zero and, "completing the square" by
adding and subtracting 
the term $\frac{\vert\nabla^\perp u\vert^2}{u}$ inside the integral, we get the formula
\begin{align*}
  \frac{d\,}{dt}\,\Bigl\{[4\pi({T}-t)]^{\frac{m-n}{2}}\int_N
    u\,d\mu_t\,\Bigr\}
 =&\,-[4\pi({T}-t)]^{\frac{m-n}{2}}\int_N
 \vert\HHH-\nabla^\perp\log{u}\vert^2u\,d\mu_t\\
&\,-[4\pi({T}-t)]^{\frac{m-n}{2}}\int_N
\Bigl(\nabla^2_{\alpha\beta}\log{u}+\frac{g_{\alpha\beta}}{2({T}-t)}\Bigr)g^{\alpha\beta}u\,d\mu_t\,,
\end{align*}
for every $t\in[0,T)$, where $\nabla^\perp$ denotes the projection on the normal space to $N$
of the gradient in $M$ of a function.

\begin{rem} 
In the special case of $M=\R^m$ and $u$ equal to the backward heat kernel
  $\rho_{x_0,T}$, the last term vanishes because of the special choice of $u$ and we 
have the "classical" Huisken's monotonicity formula.
\end{rem}

The right hand side of this formula consists of a 
nonpositive quantity (minus the integral of a perfect square times
$u$, which is positive) and a term which could be nonpositive in case the two form $\nabla^2 \log u+\frac{g}{2(T-t)}$ were nonnegative definite.

Setting $v(p,s)=u(p,T-s)$, the function $v:M\times(0,T]\to\R$
is a positive solution of the standard {\em forward} heat equation on
$(M,g)$ and setting $t={T}-s$, we have 
$\nabla^2\log{u}+\frac{g}{2({T}-t)}=
\nabla^2\log{v}+\frac{g}{2s}$. In particular, its trace (the standard
{\em full} trace) is given by 
$\Delta^M \log{v}+\frac{m}{2s}$ which is exactly the Li--Yau
quantity for positive solutions of the heat equation on a compact
manifold $(M,g)$. Actually, in the paper~\cite{liyau}, Li and Yau showed that if
the Ricci tensor of $M$ is nonnegative, then the differential Harnack 
inequality $\Delta^M \log{v}+\frac{m}{2({T}-t)}\geq 0$ holds. In the
spirit of this result, in~\cite{hamilton7} Hamilton (see
also~\cite{retomul}) generalized this  
inequality to a matrix version, showing that under the assumptions that $(M,g)$ has 
parallel Ricci tensor ($\nabla\Ric=0$) and nonnegative 
sectional curvatures, the 2--form  $\nabla^2\log{v}+g/(2s)$ is
nonnegative definite ({\em Hamilton's matrix Li--Yau Harnack
  differential inequality}).

As a consequence, under these hypotheses the two form 
$$
\nabla^2\log{u}+\frac{g}{2({T}-t)}=\nabla^2\log{v}+\frac{g}{2s}
$$
is nonnegative definite and we get Hamilton's generalization of Huisken's monotonicity formula.

\begin{teo}[Huisken's Monotonicity Formula -- Hamilton's Extension~\cite{hamilton8}]
A smooth, compact, $n$--dimensional submanifold $N$ of a Riemannian
manifold $(M,g)$ moves by mean curvature in the time interval $[0,T)$ and $u:M\times
[0,T)\to\R$ is a positive smooth solution of the backward heat
equation $u_t=-\Delta^M u$.\\
Then, if the manifold $(M,g)$ has nonnegative
sectional curvatures and satisfies $\nabla^M\Ric=0$ the quantity
$[4\pi({T}-t)]^{\frac{m-n}{2}}\int_N
    u\,d\mu_t$ is nonincreasing during the flow in $[0,T)$.

\end{teo}

\begin{rem} All this discussion in the static ambient situation
  provides a first example of the connection of the monotonicity of the ``coupled'' 
integral $[4\pi({T}-t)]^{\frac{m-n}{2}}\int_N
    u\,d\mu_t$ with the validity of a Li--Yau--Hamilton Harnack
    differential inequality.
\end{rem}

\section{Moving Ambient Spaces}

Let us now assume that the metric of the ambient space evolves according to
$\partial_tg=-2\QQQ$ (if $\QQQ=\Ric$ we have the Ricci flow) and modify the backward heat equation as follows
\begin{equation*}
u_t=-\Delta^{M}u + \KKK u
\end{equation*}
for some function $\KKK$.\\
If we repeat the previous computations in this new setting, we get two extra terms.
The first comes from the modified equation for $u$ and the
second from the effect of the motion of the ambient space on the 
time derivative of the measure $\mu_t$ induced on $N$. Indeed, the associated
metric $h_t$ on $N$ is affected not just by the motion of the 
submanifold but also by the evolution of the ambient metric $g(t)$ on $M$. 
After some computations, we have
$$
\frac{d\,}{dt}
\mu_t=(-\HHH^2-g^{ij}\QQQ_{ij})\mu_t=(-\HHH^2-{\mathrm{tr}}\,\QQQ+g^{\alpha\beta}\QQQ_{\alpha\beta})\mu_t
$$
where, as before (and in the rest of the paper), the Greek letters
$\alpha, \beta,\dots$ denote the indices associated to the 
coordinates which are normal to $N$ and 
with $i, j, k,\dots$ the indices for the coordinates on $N$.\\
With this notation, we get
\begin{align*}
\frac{d}{dt}\,\Bigl(\tau^{\frac{m-n}{2}}\int _{N}u\,d\mu_t\Bigr)
=&\,-\tau^{\frac{m-n}{2}}\int _{N}\Bigl\vert\HHH-\frac{\nabla^\perp
  u}{u}\Bigr\vert^2u\,d\mu_t\\
&\,-\tau^{\frac{m-n}{2}}
\int _{N}\Bigl(\frac{\nabla^2_{\alpha\beta}u}{u} - \frac{\nabla_{\alpha}u \nabla_{\beta}u}{u^2}
+\frac{g_{\alpha\beta}}{2\tau}\Bigr)g^{\alpha\beta}u\,d\mu_t\\
&\,+\tau^{\frac{m-n}{2}}\int _{N} (\KKK -
{\mathrm{tr}}\,\QQQ+g^{\alpha\beta}\QQQ_{\alpha\beta})\,u\,d\mu_t\\
=&\,-\tau^{\frac{m-n}{2}}\int _{N}\vert\HHH+\nabla^\perp
f\vert^2e^{-f}\,d\mu_t\\
&\,+\tau^{\frac{m-n}{2}}
\int_{N}\Bigl(\nabla^2_{\alpha\beta}f+\QQQ_{\alpha\beta}-\frac{g_{\alpha\beta}}{2\tau}\Bigr)g^{\alpha\beta}
e^{-f}\,d\mu_t\\
&\,+\tau^{\frac{m-n}{2}}\int _{N} (\KKK - {\mathrm{tr}}\,\QQQ)\,e^{-f}\,d\mu_t\,,
\end{align*}
where we substituted $\tau=T-t$ and $f=-\log{u}$, hence, $f_t=-\Delta^Mf+\vert\nabla f\vert^2-\KKK$.\\

We will concentrate on the following situations: $\QQQ=\Ric$ or
$\QQQ=-\Ric$, that is, the metric $g$ on $M$ evolves either by the Ricci flow
or by the backward Ricci flow and we will choose 
$\KKK=0$ or $\KKK=\mathrm{tr}\,\QQQ$. In this latter case the last term in the formula above clearly vanishes and we obtain
\begin{align*}
\frac{d}{dt}\Bigl(\tau^{\frac{m-n}{2}}\int _{N}u\,d\mu_t\Bigr)
=&\,-\tau^{\frac{m-n}{2}}\int _{N}\vert\HHH+\nabla^\perp
f\vert^2e^{-f}\,d\mu_t\\
&\,+\tau^{\frac{m-n}{2}}
\int_{N}\Bigl(\nabla^2_{\alpha\beta}f+\QQQ_{\alpha\beta}-\frac{g_{\alpha\beta}}{2\tau}\Bigr)
g^{\alpha\beta}e^{-f}\,d\mu_t\,.
\end{align*}
Moreover, notice that with the choice $\KKK={\mathrm{tr}}\,\QQQ$, we also have
$$
\frac{d\,}{dt}\int_M u \,= \int_M (u_t - {\mathrm{tr}}\,\QQQ u)\, =
\int_M -\Delta^Mu = 0\,,
$$
when $M$ is compact, 
hence the ``ambient'' integral $\int_M u=\int_M e^{-f}$ is constant
during the flow.

A family of metrics $g(t)$ on a manifold $M$ for $t\in[0,T)$, evolves by the {\em Ricci flow}
if $\partial_tg=-2\Ric_{g(t)}$. Moreover, we say that $g(t)$ evolves by the {\em
  backward Ricci flow} if $\partial_tg=2\Ric_{g(t)}$ (in the following the subscript $g(t)$ will be always understood).\\
 Under the Ricci flow, the Christoffel symbols of the evolving Levi--Civita connection, the Ricci tensor and the scalar curvature evolve according to
$$
\partial_t\Gamma_{ij}^k=-g^{kl}(\nabla_i\RRR_{jl}+\nabla_j\RRR_{il}-\nabla_l\RRR_{ij})\,,
$$
$$
\partial_t\RRR_{ij}=\Delta\RRR_{ij}+2\RRR^{pq}\RRR_{ipjq}-2g^{pq}\RRR_{ip}\RRR_{qj}\,,
$$
$$
\partial_t \RRR =\Delta \RRR+2\vert\Ric\vert^2\,.
$$
The analogous evolution equations for the {\em backward} Ricci flow
(inverting the time direction) are simply the same with a minus sign
in front of the right hand sides.

Here $\RRR_{ijkl}$ are the components of the $(4,0)$--Riemann tensor $\mathrm{Riem}$ (with the convention that for the standard sphere $\SS^n$ we have ${\mathrm{Riem}}(v,w,v,w)>0$), $\Ric$ is the Ricci tensor with components $\RRR_{ik}=g^{jl}\RRR_{ijkl}$ and finally $\RRR=g^{ik}\RRR_{ik}$ is the scalar curvature.

\subsection{${\mathbf {RF_0}}$ -- Ricci Flow and $\KKK=0$}

We assume $\partial_t g=-2\Ric$ and $u_t=-\Delta^Mu$, then
\begin{align*}
\frac{d}{dt}\,\Bigl(\tau^{\frac{m-n}{2}}\int _{N}u\,d\mu_t\Bigr)
=&\,-\tau^{\frac{m-n}{2}}\int _{N}\vert\HHH+\nabla^\perp
f\vert^2e^{-f}\,d\mu_t\\
&\,+\tau^{\frac{m-n}{2}}
\int_{N}\Bigl(\nabla^2_{\alpha\beta}f+\RRR_{\alpha\beta}-\frac{g_{\alpha\beta}}{2\tau}\Bigr)g^{\alpha\beta}
e^{-f}\,d\mu_t\\
&\,-\tau^{\frac{m-n}{2}}\int _{N} \RRR\,e^{-f}\,d\mu_t\\
=&\,-\tau^{\frac{m-n}{2}}\int _{N}\vert\HHH+\nabla^\perp
f\vert^2e^{-f}\,d\mu_t\\
&\,+\tau^{\frac{m-n}{2}}
\int_{N}\Bigl(\nabla^2_{\alpha\beta}f+\RRR_{\alpha\beta}-\frac{g_{\alpha\beta}}{2\tau}-\frac{\RRR
  g_{\alpha\beta}}{m-n}\Bigr)g^{\alpha\beta}e^{-f}\,d\mu_t\,,\nonumber
\end{align*}
with $f=-\log{u}$, hence, $f_t=-\Delta^Mf+\vert\nabla f\vert^2$.

\subsection{${\mathbf {BRF_0}}$ -- Back--Ricci Flow and $\KKK=0$}

We assume $\partial_t g=2\Ric$ and $u_t=-\Delta^Mu$, then
\begin{align*}
\frac{d}{dt}\,\Bigl(\tau^{\frac{m-n}{2}}\int _{N}u\,d\mu_t\Bigr)
=&\,-\tau^{\frac{m-n}{2}}\int _{N}\vert\HHH+\nabla^\perp
f\vert^2e^{-f}\,d\mu_t\\
&\,+\tau^{\frac{m-n}{2}}
\int_{N}\Bigl(\nabla^2_{\alpha\beta}f-\RRR_{\alpha\beta}-\frac{g_{\alpha\beta}}{2\tau}\Bigr)g^{\alpha\beta}e^{-f}\,d\mu_t\\
&\,+\tau^{\frac{m-n}{2}}\int _{N} \RRR\,e^{-f}\,d\mu_t\\
=&\,-\tau^{\frac{m-n}{2}}\int _{N}\vert\HHH+\nabla^\perp
f\vert^2e^{-f}\,d\mu_t\\
&\,+\tau^{\frac{m-n}{2}}
\int_{N}\Bigl(\nabla^2_{\alpha\beta}f-\RRR_{\alpha\beta}-\frac{g_{\alpha\beta}}{2\tau}+\frac{\RRR  g_{\alpha\beta}}{m-n}\Bigr)g^{\alpha\beta}e^{-f}\,d\mu_t\,,
\end{align*}
with $f=-\log{u}$, hence, $f_t=-\Delta^Mf+\vert\nabla f\vert^2$.

\subsection{${\mathbf {RF}}$ -- Ricci Flow and $\KKK={\mathrm{tr}}\,\QQQ=\RRR$}

We assume $\partial_t g=-2\Ric$ and $u_t=-\Delta^Mu+\RRR u$, then
\begin{align*}
\frac{d}{dt}\Bigl(\tau^{\frac{m-n}{2}}\int _{N}u\,d\mu_t\Bigr)
=&\,-\tau^{\frac{m-n}{2}}\int _{N}\vert\HHH+\nabla^\perp
f\vert^2e^{-f}\,d\mu_t\nonumber\\
&\,+\tau^{\frac{m-n}{2}}
\int_{N}\Bigl(\nabla^2_{\alpha\beta}f+\RRR_{\alpha\beta}-\frac{g_{\alpha\beta}}{2\tau}\Bigr)g^{\alpha\beta}
e^{-f}\,d\mu_t\,,
\end{align*}
with $f=-\log{u}$, hence $f_t=-\Delta^Mf+\vert\nabla f\vert^2-\RRR$.

Monotonicity of $\tau^{\frac{m-n}{2}}\int _{N}u\,d\mu_t$ 
is then related to the nonpositivity of the Li--Yau--Hamilton quantity 
$$
\Bigl(\nabla^2_{\alpha\beta}f
+\RRR_{\alpha\beta}-\frac{g_{\alpha\beta}}{2\tau}\Bigr)g^{\alpha\beta}\,.
$$
Notice that the same conclusion holds also if $u_t\leq-\Delta^M u + \RRR u$.

We emphasize that in the ${\mathbf {RF_0}}$ case, the same nonpositivity property 
clearly implies the monotonicity when $\RRR$ is always nonnegative.

\subsection{${\mathbf {BRF}}$ -- Back--Ricci Flow and $\KKK={\mathrm{tr}}\,\QQQ=-\RRR$}

We assume $\partial_t g=2\Ric$ and $u_t=-\Delta^Mu-\RRR u$, then
\begin{align*}
\frac{d}{dt}\Bigl(\tau^{\frac{m-n}{2}}\int _{N}u\,d\mu_t\Bigr)
=&\,-\tau^{\frac{m-n}{2}}\int _{N}\vert\HHH+\nabla^\perp
f\vert^2e^{-f}\,d\mu_t\nonumber\\
&\,+\tau^{\frac{m-n}{2}}
\int_{N}\Bigl(\nabla^2_{\alpha\beta}f-\RRR_{\alpha\beta}-\frac{g_{\alpha\beta}}{2\tau}\Bigr)g^{\alpha\beta}
e^{-f}\,d\mu_t\,,
\end{align*}
with $f=-\log{u}$, hence $f_t=-\Delta^Mf+\vert\nabla f\vert^2+\RRR$.

Monotonicity of $\tau^{\frac{m-n}{2}}\int _{N}u\,d\mu_t$ 
is then related to the nonpositivity of 
$$
\Bigl(\nabla^2_{\alpha\beta}f
-\RRR_{\alpha\beta}-\frac{g_{\alpha\beta}}{2\tau}\Bigr)g^{\alpha\beta}\,.
$$
Notice that the same conclusion holds also if $u_t\leq-\Delta^M u - \RRR u$.

\section{Ricci Solitons}\label{solsec}

We choose now $\QQQ=\Ric$, that is, the metric $g$ on $M$ evolves by the Ricci flow 
in some time  interval $I\subset\R$ 
and we set $\KKK=\RRR$ to be the scalar curvature of $(M,g)$.\\
By the previous computations in the ${\mathbf {RF}}$ case, we get
\begin{align}
\frac{d}{dt}\Bigl(\tau^{\frac{m-n}{2}}\int_{N}u\,d\mu_t\Bigr)
=&\,-\tau^{\frac{m-n}{2}}\int _{N}\vert\HHH+\nabla^\perp
f\vert^2e^{-f}\,d\mu_t\label{equ999}\\
&\,+\tau^{\frac{m-n}{2}}
\int_{N}\Bigl(\nabla^2_{\alpha\beta}f+\RRR_{\alpha\beta}-\frac{g_{\alpha\beta}}{2\tau}\Bigr)g^{\alpha\beta}
e^{-f}\,d\mu_t\,,\nonumber
\end{align}
for a positive solution of the conjugate heat equation
\begin{equation}\label{cheq}
u_t=-\Delta^M u + \RRR u
\end{equation}
and $f=-\log{u}$, $\tau=T-t$, for $t\in I$ with $t<T\in\R$. 

Let us assume that $(M,g(t))$ is a gradient soliton (a self--similar
solution) of Ricci flow and 
$F:M\times I\to\R$ its ``potential'' function, namely, 
\begin{itemize}
\item {\em Shrinking Soliton:} the flow is defined on
  $I=(-\infty,T_{\max})$, the metric $g$ and the function $F$ satisfy
  $\nabla^2 F+\Ric= g/2(T_{\max}-t)$.

\item {\em Steady Soliton:} the flow is ``eternal'', $I=\R$, 
the metric $g$ and the function $F$ satisfy $\nabla^2 F+\Ric=0$.

\item {\em Expanding Soliton:} the flow is defined on
  $I=(T_{\min},+\infty)$, 
the metric $g$ and the function $F$ satisfy 
$\nabla^2 F+\Ric=g/2(T_{\min}-t)$.
\end{itemize}

Then we analyze these three situations separately.

\begin{itemize}
\item Shrinking Solitons:
It can be seen that the function $u=e^{-F}/(T_{\max}-t)^{m/2}$ satisfies 
the conjugate heat equation~\eqref{cheq}
(see~\cite[Section~1.5]{retomul}, for instance, 
for this and the next cases). Then, letting $f=-\log{u}=F+\frac{m}{2}\log{(T_{\max}-t)}$ and 
substituting inside equation~\eqref{equ999}, we get
\begin{align*}
\frac{d}{dt}\Bigl(\tau^{\frac{m-n}{2}}\int_{N}\frac{e^{-F}}{(T_{\max}-t)^{m/2}}\,d\mu_t\Bigr)
=&\,-\frac{(T-t)^{\frac{m-n}{2}}}{(T_{\max}-t)^{m/2}}\int _{N}\vert\HHH+\nabla^\perp
F\vert^2e^{-F}\,d\mu_t\\
&\,+\frac{(T-t)^{\frac{m-n}{2}}}{(T_{\max}-t)^{m/2}}
\int_{N}\frac{m-n}{2}\Bigl(\frac{1}{T_{\max}-t}-\frac{1}{T-t}\Bigr)\,e^{-F}\,d\mu_t\,,
\end{align*}
which is nonpositive for every
$t\in(-\infty,\min\{T,T_{\max}\})$, if $T\leq T_{\max}$.\\
Actually, the right side is always negative if $T<T_{\max}$ and in
the particular case of $T=T_{\max}$, we have the neat formula
\begin{equation*}
\frac{d}{dt}\int_{N}\frac{e^{-F}}{(T_{\max}-t)^{n/2}}\,d\mu_t
=-\int _{N}\vert\HHH+\nabla^\perp
F\vert^2\frac{e^{-F}}{(T_{\max}-t)^{n/2}}\,d\mu_t\leq 0\,,
\end{equation*}
with equality if and only if the submanifold $N$ satisfies $\HHH+\nabla^\perp F=0$
at every point, for some time $t$.

An almost trivial example of this situation is a ``static'' 
maximal sphere $\SS^n$ in the sphere $\SS^m$ evolving by Ricci flow. Indeed, this latter "generates" a gradient, shrinking Ricci soliton 
with a constant in space potential function $F$ and the maximal sphere $\SS^n$ satisfies $\HHH=0$.\\
Another example is given by the flat $\R^m$ with potential function $F(x,t)=\frac{\vert x-x_0\vert^2}{4(T_{\max}-t)}$ which is called the {\em Gaussian} shrinking soliton, for some $x_0\in\R^m$. Substituting in the last equation above, one recovers the "classical" Huisken's monotonicity formula~\eqref{huiskform}.\\
Notice anyway that the family of cylinders $(\SS^2\times\R,g(t))$ with
the evolving metric $g(t)=-2t(g_{\mathrm {can}}^{\SS^2}+dr^2)$ in the
halfline $t\in(-\infty,0)$, is a gradient, shrinking, Ricci soliton with
$T_{\max}=0$ and potential function
$F:\SS^2\times\R\times(-\infty,0)\to\R$ given by
$F(\theta,r,t)=-\frac{(r-r_0)^2}{4t}$, for some $r_0\in\R$.
Any 2--sphere $\SS^2\times\{\overline{r}\}$ inside $\SS^2\times\R$ is 
actually ``static'' during its flow by mean curvature, 
since its second fundamental form (hence, its mean curvature) is zero,
but the Huisken's integral is not constant, unless $\overline{r}=r_0$ 
(it holds only for a single 2--sphere of the whole family fibering 
the cylinder). This follows easily as the vector $\nabla^\perp F=-\frac{(r-r_0)}{2t}\partial_r$ 
must be zero in such case.\\

\item Steady Solitons: The function $u=e^{-F}$ satisfies 
the conjugate heat equation~\eqref{cheq} hence, letting $f=-\log{u}=F$ in equation~\eqref{equ999} we have
$$
\frac{d}{dt}\Bigl(\tau^{\frac{m-n}{2}}\int_{N}e^{-F}\,d\mu_t\Bigr)
=-\tau^{\frac{m-n}{2}}\int _{N}\vert\HHH+\nabla^\perp
F\vert^2e^{-F}\,d\mu_t-\tau^{\frac{m-n-2}{2}}\frac{m-n}{2}
\int_{N}e^{-F}\,d\mu_t\,,
$$
which is always negative for every
$t\in(-\infty,T)$.\\
Notice that in this case, it follows
$$
\frac{d}{dt}\int_{N}e^{-F}\,d\mu_t
=-\int _{N}\vert\HHH+\nabla^\perp
F\vert^2e^{-F}\,d\mu_t\,,
$$
for every $t\in\R$.\\

\item Expanding Solitons: In this case the function 
$u=e^{-F}/(t-T_{\min})^{m/2}$ satisfies 
the conjugate heat equation~\eqref{cheq}, then, letting $f=-\log{u}=F+\frac{m}{2}\log{(t-T_{\min})}$ and 
substituting inside equation~\eqref{equ999}, we get
\begin{align*}
\frac{d}{dt}\Bigl(\tau^{\frac{m-n}{2}}\int_{N}\frac{e^{-F}}{(t-T_{\min})^{m/2}}\,d\mu_t\Bigr)
=&\,-\frac{(T-t)^{\frac{m-n}{2}}}{(t-T_{\min})^{m/2}}\int _{N}\vert\HHH+\nabla^\perp
F\vert^2e^{-F}\,d\mu_t\\
&\,+\frac{(T-t)^{\frac{m-n}{2}}}{(t-T_{\min})^{m/2}}
\int_{N}\frac{m-n}{2}\Bigl(\frac{1}{T_{\min}-t}-\frac{1}{T-t}\Bigr)\,e^{-F}\,d\mu_t\,,
\end{align*}
which is always negative for every $t\in(T_{\min},T)$ (notice that in
this case $T\leq T_{\min}$ has no meaning).
\end{itemize}

\begin{prop} If $(M,g(t))$ is an $m$--dimensional, 
shrinking, gradient Ricci soliton in the interval 
$(-\infty,T_{\max})$ and $F$ its potential function, 
then, the Huisken's integral $\tau^{\frac{m-n}{2}}\int_{N}u\,d\mu_t$,
with $u=e^{-F}/(T_{\max}-t)^{m/2}$, $\tau=T-t$ and $T\leq T_{\max}$, of an $n$--dimensional
submanifold $N$ moving by mean curvature inside $(M,g(t))$ 
is monotone nonincreasing for every $t\in(-\infty,T)$.\\
It is actually monotone decreasing, unless $T=T_{\max}$ and at some
time the submanifold $N$ satisfies $\HHH+\nabla^\perp F=0$ at every
point.

If $(M,g(t))$ is an $m$--dimensional steady or expanding, gradient
Ricci soliton with potential function $F$ in the interval
$(T_{\min},+\infty)$, then, the Huisken's integral $\tau^{\frac{m-n}{2}}\int_{N}u\,d\mu_t$, with $u=e^{-F}$ or 
$u=e^{-F}/(t-T_{\min})^{m/2}$ respectively, $\tau=T-t$, $T>T_{\min}$ and $N$ as above, 
is monotone decreasing for every $t\in(T_{\min},T)$.\\
Moreover, in the steady case, the integral $\int_{N}e^{-F}\,d\mu_t$ in
monotone nonincreasing for every $t\in\R$ and 
actually decreasing unless the submanifold $N$ satisfies $\HHH+\nabla^\perp F=0$ at every
point.
\end{prop}

\section{Computations I -- Ricci Flow and LYH Matrix Harnack Inequalities}

In this section we will deal with the ${\mathbf {RF}}$ case, that is,
we will assume that $(M,g(t))$ is an $m$--dimensional Riemannian manifold
evolving by the Ricci flow
$\partial_t g=-2\Ric$ and the smooth function $u:M\times[0,T)\to\R$ is a
positive solution of $u_t=-\Delta u+\RRR u$. Under these assumptions,
considering a compact $n$--submanifold $N$ moving by mean curvature,
we have seen that, setting $\tau=T-t$, the monotonicity of the
Huisken's integral
$$
\tau^{\frac{m-n}{2}}\int _{N}u\,d\mu_t
$$
is implied by the nonpositivity of the expression 
$$
\Bigl(\nabla^2_{\alpha\beta}f
+\RRR_{\alpha\beta}-\frac{g_{\alpha\beta}}{2\tau}\Bigr)g^{\alpha\beta}\,,
$$
with $f=-\log{u}$ which hence satisfies $f_t=-\Delta f+\vert\nabla
f\vert^2-\RRR$. This would be a straightforward consequence of the
nonpositivity (along the flow) of the full 2--form
$$
\nabla^2_{ij}f
+\RRR_{ij}-\frac{g_{ij}}{2\tau}\,,
$$
which is clearly a stronger property.\\
Equivalently, if we had chosen $f=\log u$, we would be interested in
the nonnegativity of
\begin{equation}\label{MovingHarnack}
\nabla^2_{ij}f-\RRR_{ij}+\frac{g_{ij}}{2\tau}\,,
\end{equation}
for $f=\log u$ satisfying
$$
f_t=-\Delta f - \vert\nabla f\vert^2+\RRR\,,
$$
which is an analogue of Li--Yau--Hamilton differential matrix Harnack
inequality in a moving ambient space.

We set $L_{ij}=\nabla^2_{ij}f-\RRR_{ij}$, $H_{ij}=\tau L_{ij}+
g_{ij}/2=\tau[\nabla^2_{ij}f-\RRR_{ij}]+g_{ij}/2$ and we compute the evolution
equation of the form $H$, whose nonnegativity is trivially equivalent
to the one of the form~\eqref{MovingHarnack}. In normal coordinates, using the following commutation rule between the Laplacian and
the second covariant derivatives of a function $f:M\to\R$ that can be
obtained interchanging repeatedly the covariant derivatives and using
the II Bianchi identity
\begin{align*}
\nabla^2_{ij}\Delta f-\Delta\nabla^2_{ij}f=&\,
-(\nabla_i\RRR_{jk}+\nabla_j\RRR_{ik}-\nabla_k\RRR_{ij})\nabla^k f\\
&\,-g^{pq}\RRR_{jp}\nabla^2_{qi}f-g^{pq}\RRR_{ip}\nabla^2_{qj}f
+2g^{pr}g^{qs}\RRR_{ipjq}\nabla^2_{rs}f\,,\nonumber
\end{align*}
we have
\begin{align*}
(\partial_t+\Delta)H_{ij}=&\,-L_{ij}-\RRR_{ij}\\
&\,+\tau[\Delta \nabla^2_{ij}f+\nabla^2_{ij} f_t+(\nabla_i \RRR_{jk}+\nabla_j \RRR_{ik}-\nabla_k
\RRR_{ij})\nabla_kf]\\
&\,-\tau[\partial_t\RRR_{ij}+\Delta\RRR_{ij}]\\
=&\,-L_{ij}-\RRR_{ij}\\
&\,+\tau[\Delta \nabla^2_{ij}f-\nabla^2_{ij}\Delta f-\nabla^2_{ij}\vert\nabla f\vert^2\\
&\,+(\nabla_i \RRR_{jk}+\nabla_j \RRR_{ik}-\nabla_k
\RRR_{ij})\nabla_kf]\\
&\,-\tau[2\Delta\RRR_{ij}+2\RRR_{pq}\RRR_{ipjq}-2\RRR_{ip}\RRR_{pj}-\nabla^2_{ij}\RRR]\\
=&\,-L_{ij}-\RRR_{ij}\\
&\,+\tau[(\nabla_i\RRR_{jk}+\nabla_j\RRR_{ik}-\nabla_k\RRR_{ij})\nabla_k f\\
&\,+\RRR_{jp}\nabla^2_{ip}f+\RRR_{ip}\nabla^2_{pj}f
+2\RRR_{ikpj}\nabla^2_{kp}f\\
&\,-\nabla^2_{ij}\vert\nabla f\vert^2
+(\nabla_i \RRR_{jk}+\nabla_j \RRR_{ik}-\nabla_k
\RRR_{ij})\nabla_kf]\\
&\,-\tau[2\Delta\RRR_{ij}+2\RRR_{pq}\RRR_{ipjq}-2\RRR_{ip}\RRR_{pj}-\nabla^2_{ij}\RRR]\\
=&\,-L_{ij}-\RRR_{ij}\\
&\,+\tau[\RRR_{jp}\nabla^2_{ip}f+\RRR_{ip}\nabla^2_{pj}f
+2\RRR_{ikpj}\nabla^2_{kp}f\\
&\,-2\nabla^2_{ip}f\nabla^2_{jp}f-2\nabla^3_{ijk}f\nabla_kf
+2(\nabla_i \RRR_{jk}+\nabla_j \RRR_{ik}-\nabla_k
\RRR_{ij})\nabla_kf]\\
&\,-\tau[2\Delta\RRR_{ij}+2\RRR_{pq}\RRR_{ipjq}-2\RRR_{ip}\RRR_{pj}-\nabla^2_{ij}\RRR]\,.
\end{align*}
Commuting the covariant derivatives of the term containing the third derivatives of $f$, that is, $\nabla^3_{ijk}f=\nabla^3_{kij}f+\RRR_{ikjp}\nabla_p f$, we get
\begin{align*}
(\partial_t+\Delta)H_{ij}=&\,-L_{ij}-\RRR_{ij}\\
&\,+\tau[\RRR_{jp}\nabla^2_{ip}f+\RRR_{ip}\nabla^2_{pj}f-2\nabla^2_{ip}f\nabla^2_{jp}f-2\nabla^3_{kij}f\nabla_kf]\\
&\,-\tau[2\Delta\RRR_{ij}+2\RRR_{pq}\RRR_{ipjq}-2\RRR_{ip}\RRR_{pj}-\nabla^2_{ij}\RRR]\\
&\,+\tau[2(\nabla_i \RRR_{jk}+\nabla_j \RRR_{ik}-\nabla_k
\RRR_{ij})\nabla_kf\\
&\,-2\RRR_{ikjp}\nabla^2_{pk}f-2\RRR_{ikjp}\nabla_p f\nabla_k f]\,.
\end{align*}
Finally, substituting $L_{ij}=[H_{ij}-g_{ij}/2]/\tau$ and $\nabla^2_{ij}f=[H_{ij}-g_{ij}/2]/\tau+\RRR_{ij}$, 
we obtain
\begin{align*}
(\partial_t+\Delta)H_{ij}
=&\,[H_{ij}-2H^2_{ij}]/\tau-2\nabla_k H_{ij}\nabla_kf\\
&\,-[\RRR_{ik}H_{jk}+\RRR_{jk}H_{ik}+2\RRR_{ikjp}H_{pk}]\\
&\,-\tau[2\Delta\RRR_{ij}-2\RRR_{jp}\RRR_{ip}+4\RRR_{pq}\RRR_{ipjq}-\nabla^2_{ij}\RRR-\RRR_{ij}/\tau]\\
&\,+\tau[2(\nabla_i \RRR_{jk}+\nabla_j \RRR_{ik}-2\nabla_k\RRR_{ij})\nabla_kf]\\
&\,-2\tau\RRR_{ikjp}\nabla_p f\nabla_k f\,,
\end{align*}
which in a generic coordinate system reads
\begin{align*}
(\partial_t+\Delta)H_{ij}
=&\,[H_{ij}-2H^2_{ij}]/\tau-2\nabla_k H_{ij}\nabla^kf-g^{pq}\RRR_{ip}H_{jq}-g^{pq}\RRR_{jp}H_{iq}-2\RRR_{ipjq}H^{pq}\\
&\,-\tau\Bigl[2\Delta\RRR_{ij}-\nabla^2_{ij}\RRR-2g^{pq}\RRR_{ip}\RRR_{jq}+4\RRR^{pq}\RRR_{ipjq}-\RRR_{ij}/\tau\phantom{\Bigl]}\\
&\,\phantom{-\tau\Bigl[ \,}-2(\nabla_i \RRR_{jk}+\nabla_j \RRR_{ik}-2\nabla_k\RRR_{ij})\nabla^kf +2\RRR_{ipjq}\nabla^p f\nabla^q f\Bigr]\\
=&\,[H_{ij}-2H^2_{ij}]/\tau-2\nabla_k H_{ij}\nabla^kf-\RRR_i^kH_{kj}-\RRR_j^kH_{ki}-2\RRR_{ipjq}H^{pq}-\tau W_{ij}\,,
\end{align*}
where we set
\begin{align*}
W_{ij}=&\,2\Delta\RRR_{ij}-\nabla^2_{ij}\RRR-2g^{pq}\RRR_{ip}\RRR_{jq}+4\RRR^{pq}\RRR_{ipjq}-\RRR_{ij}/\tau\\
&\,-2(\nabla_i \RRR_{jk}+\nabla_j \RRR_{ik}-2\nabla_k\RRR_{ij})\nabla^kf +2\RRR_{ipjq}\nabla^p f\nabla^q f\,.
\end{align*}
Notice that when $t>T$ (not our case) this form $W$ is the Hamilton's Harnack quadratic, 
defined in~\cite{hamilton12}, contracted with $\nabla f$, 
(this term with the ``wrong time'' also appears in the computations
about the {\em reduced length} in Perelman's paper~\cite{perel1}). 
This quantity vanishes on a shrinking, gradient Ricci
soliton with $T=T_{\max}$ when $f$ is equal to minus its potential function $F$, 
so sometimes it is called Hamilton's matrix Harnack quadratic 
{\em for shrinkers} (the original Hamilton's Harnack quadratic is
instead zero on expanders).

Arguing as in Hamilton~\cite{hamilton7} by means of his matrix
maximum principle, if there is a sequence $t_i\to T$ such that the form
$H$ is positive definite, the Riemann curvature operator and the
form $W$ are nonnegative definite in $M\times[0,T)$, then it follows that the form $H$ is
nonnegative definite in the whole $M\times[0,T)$ which is what we 
need to get the monotonicity of the Huisken's integral.\\ 
Unfortunately, the form $W$ is not, in general, nonnegative definite for every function $f$, 
even if the Ricci flow is {\em ancient} and the Riemann curvature
operator is nonnegative, in contrast to the nicely behaved original
Hamilton's Harnack quadratic. Anyway, when the flow is a gradient,
shrinking Ricci soliton with nonnegative Riemann operator and
$\tau=T_{\max}-t$, 
the form $W$ is nonnegative definite for every function $f$, 
indeed, there hold (by the soliton equation, 
see~\cite[Chapter~8, Section~5]{chowluni})
$$
2\Delta\RRR_{ij}-\nabla^2_{ij}\RRR-2g^{pq}\RRR_{ip}\RRR_{jq}+4\RRR^{pq}\RRR_{ipjq}-\RRR_{ij}/\tau
=2(\nabla_k\RRR_{ij}-2\nabla_j\RRR_{ik})\nabla^kF
$$
and
$$
\nabla_j\RRR_{ki}-\nabla_k\RRR_{ji}=-\RRR_{jkip}\nabla^pF\,,
$$
hence, the equality
\begin{align*}
W_{ij}=&\,2\RRR_{jkip}\nabla^pF\nabla^k F
+2(\RRR_{jkip}\nabla^pF+\RRR_{ikjp}\nabla^pF)\nabla^kf +2\RRR_{ikjp}\nabla^p f\nabla^k f\\
=&\,2\RRR_{jkip}\nabla^p(F+f)\nabla^k(F+f)\,,
\end{align*}
implies the claim, by the curvature assumption. Then, 
if a solution $u$ of the conjugate heat equation~\eqref{cheq} satisfies 
$$
\nabla^2\log{u(\cdot,t_i)}-\Ric(\cdot,t_i)+\frac{g(\cdot,t_i)}{2(T_{\max}-t)}\geq 0\,,
$$
for a sequence of times $t_i\to T_{\max}$ on the whole $M$, the monotonicity of the Huisken's
integral follows in the interval $[0,T_{\max})$.\\
With an analogous argument, it can be shown that if the flow is a gradient, steady Ricci soliton with
nonnegative Riemann operator and the function $u$ satisfies the same condition as before, then the
monotonicity of the Huisken's integral holds in the interval $[0,T)$ (in this
case the interval can also be $\R$ and $T=+\infty$).

\section{Computations II -- Backward Ricci Flow}

We deal now with the ${\mathbf {BRF}}$ case, that is, $(M,g(t))$ is an
$m$--dimensional Riemannian manifold evolving by the backward Ricci flow 
$\partial_t g=2\Ric$ and the smooth function $u:M\times[0,T)\to\R$ is a
positive solution of $u_t=-\Delta u-\RRR u$. Then, we have seen that if
$\tau=T-t$ the monotonicity of the Huisken's integral
$$
\tau^{\frac{m-n}{2}}\int _{N}u\,d\mu_t\,,
$$
where $N$ is a compact $n$--submanifold moving by mean curvature, 
is implied by the nonpositivity of the expression 
$$
\Bigl(\nabla^2_{\alpha\beta}f
-\RRR_{\alpha\beta}-\frac{g_{\alpha\beta}}{2\tau}\Bigr)g^{\alpha\beta}\,,
$$
with $f=-\log{u}$ which hence satisfies $f_t=-\Delta f+\vert\nabla
f\vert^2+\RRR$. Choosing instead $f=\log u$ which then satisfies
$$
f_t=-\Delta f - \vert\nabla f\vert^2-\RRR
$$
the above monotonicity would be a consequence of the stronger
statement that the full 2--form
$$
\nabla^2_{ij}f+\RRR_{ij}+\frac{g_{ij}}{2\tau}
$$
is nonnegative definite.

We then set $L_{ij}=\nabla^2_{ij}f+\RRR_{ij}$, $H_{ij}=\tau L_{ij}+
g_{ij}/2=\tau[\nabla^2_{ij}f+\RRR_{ij}]+g_{ij}/2$ and we compute the evolution
equation of the form $H$ (as before) whose nonnegativity is trivially equivalent
to the one of the form above. In normal coordinates, we have (along
the same line of the ${\mathbf {RF}}$ case)
\begin{align*}
(\partial_t+\Delta)H_{ij}=&\,-L_{ij}+\RRR_{ij}\\
&\,+\tau[\Delta \nabla^2_{ij}f+\nabla^2_{ij} f_t-(\nabla_i \RRR_{jk}+\nabla_j \RRR_{ik}-\nabla_k
\RRR_{ij})\nabla_kf]\\
&\,+\tau[\partial_t\RRR_{ij}+\Delta\RRR_{ij}]\\
=&\,-\nabla^2_{ij}f+\tau[\Delta \nabla^2_{ij}f-\nabla^2_{ij}\Delta f-\nabla^2_{ij}\vert\nabla
f\vert^2-(\nabla_i \RRR_{jk}+\nabla_j \RRR_{ik}-\nabla_k
\RRR_{ij})\nabla_kf]\\
&\,-\tau[2\RRR_{pq}\RRR_{ipjq}-2\RRR_{ip}\RRR_{pj}+\nabla^2_{ij}\RRR]\\
=&\,-\nabla^2_{ij}f+\tau[
(\nabla_i\RRR_{jk}+\nabla_j\RRR_{ik}-\nabla_k\RRR_{ij})\nabla_k f
+\RRR_{jp}\nabla^2_{ip}f+\RRR_{ip}\nabla^2_{pj}f
-2\RRR_{ipjq}\nabla^2_{pq}f]\\
&\,+\tau[-\nabla^2_{ij}\vert\nabla f\vert^2-(\nabla_i \RRR_{jk}+\nabla_j \RRR_{ik}-\nabla_k
\RRR_{ij})\nabla_kf]\\
&\,-\tau[2\RRR_{pq}\RRR_{ipjq}-2\RRR_{ip}\RRR_{pj}+\nabla^2_{ij}\RRR]\\
=&\,-\nabla^2_{ij}f+\tau[\RRR_{jp}\nabla^2_{ip}f+\RRR_{ip}\nabla^2_{pj}f
-2\RRR_{ipjq}\nabla^2_{pq}f-\nabla^2_{ij}\vert\nabla f\vert^2]\\
&\,-\tau[2\RRR_{pq}\RRR_{ipjq}-2\RRR_{ip}\RRR_{pj}+\nabla^2_{ij}\RRR]\\
=&\,-\nabla^2_{ij}f+\tau[\RRR_{jp}\nabla^2_{ip}f+\RRR_{ip}\nabla^2_{pj}f
-2\RRR_{ipjq}\nabla^2_{pq}f]\\
&\,-\tau[2\RRR_{pq}\RRR_{ipjq}-2\RRR_{ip}\RRR_{pj}+\nabla^2_{ij}\RRR]\\
&\,-\tau[2\nabla^2_{ip}f\nabla^2_{jp}f+2\nabla^3_{ijk}f\nabla_kf]\\
=&\,-\nabla^2_{ij}f+\tau[\RRR_{jp}\nabla^2_{ip}f+\RRR_{ip}\nabla^2_{pj}f
-2\RRR_{ipjq}\nabla^2_{pq}f]\\
&\,-\tau[2\RRR_{pq}\RRR_{ipjq}-2\RRR_{ip}\RRR_{pj}+\nabla^2_{ij}\RRR]\\
&\,-\tau[2\nabla^2_{ip}f\nabla^2_{jp}f+2\nabla^3_{kij}f\nabla_kf
+2\RRR_{ipjq}\nabla_pf\nabla_qf]\,.
\end{align*}
Substituting now $L_{ij}=[H_{ij}-g_{ij}/2]/\tau$ and $\nabla^2_{ij}f=[H_{ij}-g_{ij}/2]/\tau-\RRR_{ij}$, 
we get
\begin{align*}
(\partial_t+\Delta)H_{ij}
=&\,-H_{ij}/\tau+g_{ij}/2\tau+\RRR_{ij}\\
&\,-\tau[\RRR_{jp}\RRR_{ip}+\RRR_{ip}\RRR_{pj}-2\RRR_{ipjq}\RRR_{pq}]\\
&\,+ [\RRR_{jp}H_{ip}+\RRR_{ip}H_{pj}-2\RRR_{ipjq}H_{pq}]\\
&\,-\tau[2\RRR_{pq}\RRR_{ipjq}-2\RRR_{ip}\RRR_{pj}+\nabla^2_{ij}\RRR]\\
&\,-2\tau[H^2_{ij}/\tau^2-H_{ij}/\tau^2+g_{ij}/4\tau^2+\RRR_{ik}\RRR_{kj}
-\RRR_{ik}H_{jk}/\tau-\RRR_{jk}H_{ik}/\tau+\RRR_{ij}/\tau]\\
&\,-2\tau[\nabla_kH_{ij}\nabla_kf/\tau-\nabla_k\RRR_{ij}\nabla_kf
+\RRR_{ipjq}\nabla_pf\nabla_qf]\\
=&\,[H_{ij}-2H^2_{ij}]/\tau-\RRR_{ij}\\
&\,+ [3\RRR_{jp}H_{ip}+3\RRR_{ip}H_{pj}-2\RRR_{ipjq}H_{pq}]\\
&\,-\tau\nabla^2_{ij}\RRR-2\tau\RRR_{ik}\RRR_{kj}\\
&\,-2\tau[\nabla_kH_{ij}\nabla_kf/\tau-\nabla_k\RRR_{ij}\nabla_kf
+\RRR_{ipjq}\nabla_pf\nabla_qf]\\
=&\,[H_{ij}-2H^2_{ij}]/\tau-2\nabla_kH_{ij}\nabla_kf+3\RRR_{jp}H_{ip}+3\RRR_{ip}H_{pj}-2\RRR_{ipjq}H_{pq}\\
&\,-\tau[\nabla^2_{ij}\RRR+2\RRR_{ik}\RRR_{kj}+\RRR_{ij}/\tau-2\nabla_k\RRR_{ij}\nabla_kf+2\RRR_{ipjq}\nabla_pf\nabla_qf]\,.
\end{align*}
Thus, getting back to generic coordinates
\begin{align*}
(\partial_t+\Delta)H_{ij}
=&\,[H_{ij}-2H^2_{ij}]/\tau-2\nabla_kH_{ij}\nabla^kf+3g^{pq}\RRR_{ip}H_{jq}+3g^{pq}\RRR_{jp}H_{iq}-2\RRR_{ipjq}H^{pq}\\
&\,-\tau[\nabla^2_{ij}\RRR+2g^{pq}\RRR_{ip}\RRR_{jq}+\RRR_{ij}/\tau-2\nabla_k\RRR_{ij}\nabla^kf+2\RRR_{ipjq}\nabla^pf\nabla^qf]\\
=&\,[H_{ij}-2H^2_{ij}]/\tau-2\nabla_kH_{ij}\nabla^kf+3g^{pq}\RRR_{ip}H_{jq}+3g^{pq}\RRR_{jp}H_{iq}-2\RRR_{ipjq}H^{pq}\\
&\,-\tau Z_{ij}\,,
\end{align*}
where we set
$$
Z_{ij}=\nabla^2_{ij}\RRR+2g^{pq}\RRR_{ip}\RRR_{jq}+\RRR_{ij}/\tau-2\nabla_k\RRR_{ij}\nabla^kf+2\RRR_{ipjq}\nabla^pf\nabla^qf\,.
$$
Then, arguing now as in the  ${\mathbf {RF}}$ case, 
assuming that the Riemann curvature operator is nonnegative
(such a condition is not preserved in general under the backward Ricci
flow), if the form $Z$ is nonnegative definite we can conclude that 
if there is a sequence $t_i\to T$ such that the form
$H(\cdot,t_i)$ is positive definite, then the form $H$ is nonnegative definite on the whole
$M\times(0,T]$ and the monotonicity of the Huisken's integral follows.\\
Notice that the trace of $Z_{ij}$,
$$
g^{ij}Z_{ij}=\Delta\RRR+2\vert\Ric\vert^2 +\RRR/\tau-2\nabla_k\RRR\nabla^kf+2\RRR_{pq}\nabla^pf\nabla^qf
$$
coincides with the trace of the original Hamilton's Harnack quadratic
\begin{align*}
2&\,\Delta\RRR_{ij}-\nabla^2_{ij}\RRR-2g^{pq}\RRR_{ip}\RRR_{jq}+4\RRR^{pq}\RRR_{ipjq}+\RRR_{ij}/\tau\\
&\,-2(\nabla_i \RRR_{jk}+\nabla_j \RRR_{ik}-2\nabla_k\RRR_{ij})\nabla^kf +2\RRR_{ipjq}\nabla^p f\nabla^q f\,,
\end{align*}
after changing the sign of the function $f$. Hence, one can ask himself if under the backward Ricci flow of a
manifold with nonnegative definite Riemann curvature operator, the 2--forms
$$
Z_{ij}^U=\nabla^2_{ij}\RRR+2\RRR_{ij}^2+\RRR_{ij}/\tau-2\nabla_k\RRR_{ij} U^k
+2\RRR_{ipjq}U^pU^q
$$
are all nonnegative definite, for every vector $U=\{U^i\}$ (see Ni~\cite[Remark~6.4]{leini4}).\\
Unfortunately, this does not hold even in
dimension two, indeed, in such case we have $\RRR_{ij}=\RRR g_{ij}/2$
and $\RRR_{ijkl}=\RRR(g_{ik}g_{jl}-g_{il}g_{jk})/2$, hence, the
expression for $Z^U_{ij}$ becomes
$$
Z^U_{ij}=\nabla^2_{ij}\RRR+\frac{\RRR^2}{2}g_{ij}+\frac{\RRR}{2\tau}g_{ij}
-\langle\nabla\RRR\,\vert\,U\rangle g_{ij}+\RRR\vert
U\vert^2g_{ij}-\RRR U_i U_j\,.
$$
Checking the 2--form $Z^{\widetilde{U}}_{ij}$, where
$\widetilde{U}=\lambda U$ for $\lambda\in\R$, against the
vector $U$ we get
$$
Z^{\widetilde{U}}_{ij}U^iU^j=\Bigl[\nabla^2_{ij}\RRR U^i
U^j+\Bigl(\frac{\RRR^2}{2}+\frac{\RRR}{2\tau}\Bigr)\vert
U\vert^2\Bigr]
-\lambda\langle\nabla\RRR\,\vert\,U\rangle\vert U\vert^2\,.
$$
Therefore, if $\RRR$ is not constant, choosing $U=\nabla\RRR$, when
$\lambda>0$ is large enough this expression is negative somewhere.

\subsection{A Very Special Case} 

In dimension 2, for a surface with
positive scalar curvature, the function $u=\RRR>0$ satisfies 
$$
u_t=-\Delta u-\RRR u\,.
$$
Indeed, under the backward Ricci flow, we have
$$
\partial_t \RRR =-\Delta \RRR-\RRR^2\,,
$$
hence, the scalar curvature is a solution of the conjugate 
heat equation in dimension two (under the backward Ricci flow).\\
Then, for a closed curve $\gamma$ evolving by its curvature inside a
surface moving by backward Ricci flow, we have  
$$
\frac{d}{dt}\left(\sqrt{\tau}\int_{\gamma}\RRR\,d\mu_t\right)
=-\sqrt{\tau}\int_{\gamma}\left\vert\HHH-\nabla^\perp\log\RRR\right\vert^2\RRR\,d\mu_t
-\sqrt{\tau}\int_{\gamma} \left(\nabla^2_{\nu\nu}\log\RRR+\frac{\RRR}{2}+\frac{1}{2\tau}\right)\RRR
\,d\mu_t\,,
$$
where $\nu$ is the unit normal to the curve.\\
In this situation, the Li--Yau quadratic
$$
\nabla^2_{\nu\nu}\log\RRR+\frac{\RRR}{2}+\frac{1}{2\tau}
$$
is nonnegative, being exactly the ``special'' form of
Hamilton's Harnack inequality for surfaces with bounded positive scalar
curvature (see~\cite[Proposition~15.10]{chowbookII}) evaluated on the
pair of vectors $(\nu,\nu)$.

\begin{prop}
If $(M,g(t))$ is a family of surfaces with bounded positive scalar curvature
$\RRR$ moving by backward Ricci flow and $\gamma$ is a curve moving by
its curvature inside $(M,g(t))$, we have
$$
\frac{d}{dt}\left(\sqrt{\tau}\int_{\gamma}\RRR\,d\mu_t\right)
\leq -\sqrt{\tau}\int_{\gamma}\left\vert\HHH-\nabla^\perp\log\RRR\right\vert^2\RRR\,d\mu_t\,.
$$
\end{prop}

The inequality becomes an equality if and only if $M$ is a 
gradient, expanding Ricci soliton with $\RRR>0$ and
$\kkk=\nabla^\perp\log\RRR$ (see~\cite[Chapter~15]{chowbookII}).

\bigskip

\begin{ackn} 
The first author was partially supported by the German research
group ``Forschergruppe DFG 718, Analysis and stochastics in complex
physical systems''. 
The second author is partially supported by the Italian project
FIRB--IDEAS ``Analysis and Beyond''.\\
Some results in this paper were presented in the note~\cite{stratos1} by the third author.\\
We thank the anonymous referees for several corrections and valuable suggestions.
\end{ackn}

\bibliographystyle{amsplain}
\bibliography{biblio}

\providecommand{\bysame}{\leavevmode\hbox to3em{\hrulefill}\thinspace}
\providecommand{\MR}{\relax\ifhmode\unskip\space\fi MR }
\providecommand{\MRhref}[2]{%
  \href{http://www.ams.org/mathscinet-getitem?mr=#1}{#2}
}
\providecommand{\href}[2]{#2}
\begin{thebibliography}{10}

\bibitem{chowbookII}
B.~Chow, S.-C. Chu, D.~Glickenstein, C.~Guenther, J.~Isenberg, T.~Ivey,
  D.~Knopf, P.~Lu, F.~Luo, and L.~Ni, \emph{The {R}icci flow: techniques and
  applications. {P}art {II}. {A}nalytic aspects}, Mathematical Surveys and
  Monographs, vol. 144, American Mathematical Society, Providence, RI, 2008.

\bibitem{chowluni}
B.~Chow, P.~Lu, and L.~Ni, \emph{Hamilton's {R}icci flow}, Graduate Studies in
  Mathematics, vol.~77, American Mathematical Society, Providence, RI, 2006.

\bibitem{eck3}
K.~Ecker, \emph{A formula relating entropy monotonicity to {H}arnack
  inequalities}, Comm. Anal. Geom. \textbf{15} (2007), no.~5, 1025--1061.

\bibitem{hamilton12}
R.~S. Hamilton, \emph{The {H}arnack estimate for the {R}icci flow}, J. Diff.
  Geom. \textbf{37} (1993), no.~1, 225--243.

\bibitem{hamilton7}
\bysame, \emph{A matrix {H}arnack estimate for the heat equation}, Comm. Anal.
  Geom. \textbf{1} (1993), no.~1, 113--126.

\bibitem{hamilton8}
\bysame, \emph{Monotonicity formulas for parabolic flows on manifolds}, Comm.
  Anal. Geom. \textbf{1} (1993), no.~1, 127--137.

\bibitem{huisk3}
G.~Huisken, \emph{Asymptotic behavior for singularities of the mean curvature
  flow}, J. Diff. Geom. \textbf{31} (1990), 285--299.

\bibitem{liyau}
P.~Li and S.-T. Yau, \emph{On the parabolic kernel of the {S}chr\"odinger
  operator}, Acta Math. \textbf{156} (1986), no.~3--4, 153--201.

\bibitem{lott2}
J.~Lott, \emph{Mean curvature flow in a {R}icci flow background}, Comm. Math.
  Phys. \textbf{313} (2012), no.~2, 517--533.

\bibitem{retomul}
R.~M{\"u}ller, \emph{Differential {H}arnack inequalities and the {R}icci flow},
  EMS Series of Lectures in Mathematics, European Mathematical Society (EMS),
  Z\"urich, 2006.

\bibitem{leini4}
L.~Ni, \emph{A matrix {L}i--{Y}au--{H}amilton estimate for {K}\"ahler-{R}icci
  flow}, J. Diff. Geom. \textbf{75} (2007), no.~2, 303--358.

\bibitem{perel1}
G.~Perelman, \emph{The entropy formula for the {R}icci flow and its geometric
  applications}, ArXiv Preprint Server -- http://arxiv.org, 2002.

\bibitem{stratos1}
E.~Tsatis, \emph{Mean curvature flow on {R}icci solitons}, J. Phys. A
  \textbf{43} (2010), no.~4, 045202, 13.

\end{thebibliography}

\end{document}